\documentclass[12pt,a4paper
]{article}
\usepackage[T2A]{fontenc}
\usepackage[cp1251]{inputenc}
\usepackage[russian]{babel}
\usepackage{srcltx}
\usepackage{mathtext}
\usepackage{cmap}
\usepackage{latexsym,amsfonts,amssymb,amsmath,longtable,
mathrsfs,theorem,hyperref
}
\usepackage{mathdots}
\usepackage{graphicx}
\usepackage[pic]{xy}
\renewcommand{\le}{\leq}

\renewcommand{\mathcal}{\mathscr}
\theoremstyle{plain}
{
\newtheorem{Probl}{{\bfseries Проблема}}

{\theorembodyfont{\itshape}
\newtheorem{Lemma}{{\bfseries Лемма}}}
{\theorembodyfont{\itshape}
\newtheorem{Cor}{{\bfseries Следствие}}
\newtheorem{Theo}{{\bfseries Теорема}}
\newtheorem{Prop}{{\bfseries Предложение}}

}

\DeclareMathOperator{\m}{m} 
\DeclareMathOperator{\Irr}{Irr}  \DeclareMathOperator{\Cay}{Cay}

\setlength{\textwidth}{160mm}
\setlength{\textheight}{250mm} \headheight0mm \headsep0mm

\newcommand{\X}{\mathfrak X}

\title{\vspace{-1cm} \hfill{\normalsize  }{
\fontfamily{cmr} \fontseries{bx} \selectfont \\
\vspace{1cm} 
Spectra of Cayley graphs}
\thanks{Исследования первого автора поддержаны Национальным естественно-научным фондом Китая (NNSF of China), грант \#11771409.
Исследования четвертого автора  поддержаны Президентом Китайской академии наук (CAS President’s International Fellowship Initiative, PIFI), грант \#2016VMA078.}}
\date{}
\author{\bf   
W.~Guo, D.~V.~Lytkina, V.~D.~Mazurov, D.~O.~Revin}

\begin{document}

\newcommand{\qed}{\hfill{$\square$}}
\renewcommand{\proofname}{{ {\rm ДОКАЗАТЕЛЬСТВО}}}

%
%
%


\maketitle
\pagenumbering{arabic}
\begin{center}
{\bf {\small Abstract}}
\end{center}

{\small

Let $G$~be a group and $S\subseteq G$ its subset such that $S=S^{-1}$, where $S^{-1}=\{s^{-1}\mid s\in S\}$. Then {\it the Cayley graph ${\rm Cay}(G,S)$} is an undirected graph $\Gamma$
with the vertex set $V(\Gamma)=G$ and the edge set $E(\Gamma)=\{(g,gs)\mid g\in G, s\in S\}$. 
A graph $\Gamma$ is said to be {\it integral} if every eigenvalue of the adjacency matrix of
$\Gamma$ is integer.
In the paper, we prove the following theorem: { \it if a subset $S=S^{-1}$ of $G$ is  normal and  $s\in S\Rightarrow s^k\in S$ for every $k\in \mathbb{Z}$ such that $(k,|s|)=1$,
then
${\rm Cay}(G,S)$ is integral.} In particular, {\it if $S\subseteq G$ is a normal set of involutions, then  ${\rm Cay}(G,S)$ is integral.}
 We also use the theorem to prove that {\it if $G=A_n$ and $S=\{(12i)^{\pm1}\mid i=3,\dots,n\}$,  then
${\rm Cay}(G,S)$ is integral.} Thus, we give positive solutions for both problems 19.50(a) and 19.50(b) in ``Kourovka Notebook''.
 \medskip

}

\section{Введение}
Мы рассматриваем только конечные группы и  графы и используем термины ``группа'' и ``граф'' в значении ``конечная группа'' и ``конечный граф''.
Символами $G$ и $\Gamma$ соответственно всегда обозначаются некоторая группа  и некоторый граф. При этом $V=V(\Gamma)$ и $E=E(\Gamma)$ --- множества вершин и ребер графа~$\Gamma$.
{\it Матрицей смежности} графа $\Gamma$ с множеством вершин $V$ и множеством ребер
$E$ называется матрица $(a_{ij})\in {\rm M}_{|V|}(\mathbb{C})$, в которой строки и столбцы индексированы элементами множества $V$ и
$$
a_{ij}=\left\{\begin{array}{cc}
                1, & (i,j)\in E; \\
                0, & (i,j)\notin E
              \end{array}
\right.\quad \text{ для всех }\quad i,j\in V.
$$
Граф $\Gamma$ называется {\it целочисленным}, если его {\it спектр}, т.~е. спектр его матрицы смежности, состоит из целых чисел.

Пусть $S\subseteq G$. {\it Граф Кэли ${\rm Cay}(G,S)$ группы $G$, ассоциированный с множеством}~$S$,~--- это граф $\Gamma$
с множеством вершин $V(\Gamma)=G$ и множеством ребер $E(\Gamma)=\{(g,gs)\mid g\in G, s\in S\}$. Далее мы будем предполагать, что множество $S$, с которым ассоциирован граф ${\rm Cay}(G,S)$,
\begin{itemize}
 \item не содержит 1, что
равносильно отсутствию петель в ${\rm Cay}(G,S)$, и
\item {\it симметрично}, т.~е. $S=S^{-1}$, где $S^{-1}=\{s^{-1}\mid s\in S\}$.
\end{itemize}
Последнее условие позволяет нам рассматривать граф ${\rm Cay}(G,S)$ как неориентированный граф.

 Мы будем рассматривать только неориентированные графы без петель. Матрица смежности такого графа симметрична, и поэтому ее спектр состоит из действительных чисел.

Нам понадобятся еще два определения для формулировки основных результатов.

Подмножество $S\subseteq G$ называется {\it нормальным}, если $S=S^G=\{s^g\mid g\in G\}$ (здесь $s^g=g^{-1}sg$ --- элемент, сопряженный с $s$ с помощью элемента $g\in G$).

Известно, что для элемента $s\in G$ порядка $m$ множество порождающих элементов циклической группы $\langle s\rangle$ совпадает с множеством $$\{s^k\mid k\in \mathbb{Z}, (k,m)=1\}=\{s^k\mid 0\le k\le m-1, (k,m)=1\}$$ и его
мощность равна $\phi(m)$, где $\phi$~--- функция Эйлера, принимающая на натуральном числе $m$ с каноническим примарным разложением $m=p_1^{\alpha_1}\dots p_t^{\alpha_t}$ значение
$$
\phi(m)=p_1^{\alpha_1-1}(p_1-1)\dots p_t^{\alpha_t-1}(p_t-1).
$$
Подмножество $S\subseteq G$ назовем {\it эйлеровым}, если множество $$\{x\in G\mid \langle x\rangle=\langle s\rangle\}=\{s^k\mid 0\le k\le |s|-1, (k,|s|)=1\}$$ содержится в  $S$  для любого $s\in S.$
Заметим, что всякое эйлерово множество симметрично, поскольку $\langle s\rangle=\langle s^{-1}\rangle$ для любого $s\in G$. Если симметричное множество $S \subseteq G$ таково, что $\phi(|s|)\le 2$ или,
что то же самое,
$|s|\in\{2,3,4,6\}$ для любого $s\in S$, то оно эйлерово.

В ``Коуровской тетради'' \cite{Kour} записаны следующие две проблемы.

\begin{Probl}\label{NormInvolutions} {\rm \cite[вопрос~19.50(a)]{Kour}}
 Верно ли, что если $S\subseteq G$~--- нормальное множество элементов порядка $2$, то граф ${\rm Cay}(G,S)$ целочисленный?
\end{Probl}

\begin{Probl}\label{3Cycles} {\rm  \cite[вопрос~19.50(b)]{Kour}}
 Верно ли, что если $A_n$~--- знакопеременная группа степени~$n$ и $S=\{(12i)^{\pm1}\mid i=3,\dots,n\}$,  то граф ${\rm Cay}(A_n,S)$ целочисленный?
\end{Probl}

В данной заметке мы выведем положительный ответ на оба эти вопроса из следующего утверждения.

\begin{Theo}\label{main}  {\it Если $S$~--- эйлерово нормальное подмножество в группе $G$, то граф ${\rm Cay}(G,S)$ целочисленный.}
  \end{Theo}

  Доказательство этой теоремы дается методами теории характеров и не зависит от каких-либо других результатов о целочисленности графов.

  Частным случаем теоремы~\ref{main} оказывается \cite[теорема~1]{KonLyt}~--- результат, аналогичный теореме~\ref{main}, но c дополнительным предположением о нильпотентности группы~$G$.
  \medskip


  Положительное решение проблемы~\ref{NormInvolutions} дает
  \begin{Cor}\label{Cor1} {\it Если порядки элементов нормального симметричного  подмножества $S$ в группе $G$  принадлежат множеству
$\{2,3,4,6\}$, то граф ${\rm Cay}(G,S)$ целочисленный.}
  \end{Cor}

  Частным случаем нормального множества элементов порядка 2 (которое автоматически будет симметричным) является множество всех транспозиций в симметрической группе. Поэтому мы получаем независимое доказательство
  следующего утверждения, установленного в \cite[теорема~2]{KonLyt} с помощью целочисленности так называемого {\it звездного графа} ${\rm Cay}(G,S)$, где $G=S_n$~--- симметрическая группа степени~$n$, а
  $S=\{(1i)\mid 1< i\le n\}$ \cite[теорема~1]{KrakMoh}, \cite[следствие~2.1]{ChapFer}.

\begin{Cor}\label{Cor2} {\it Граф ${\rm Cay}(G,S)$, где $G=S_n$, а
  $S=\{(ij)\mid 1\le i< j\le n\}$, является целочисленным.}
  \end{Cor}

  Независимое от \cite[теорема~1]{KrakMoh} и \cite[следствие~2.1]{ChapFer} доказательство целочисленности самого звездного графа вытекает из следствия~\ref{Cor2} и еще одного следствия теоремы~\ref{main}:
 \begin{Cor}\label{Cor3} {\it  Пусть $R$~--- эйлерово нормальное подмножество в группе $G$ и $H$~--- подгруппа в $G$. Положим $S=R\setminus(R\cap H)$. Тогда граф ${\rm Cay}(G,S)$ целочисленный.}
  \end{Cor}

\begin{Cor}\label{Cor4} {\it Граф ${\rm Cay}(G,S)$, где $G=S_n$, а
  $S=\{(1i)\mid 1< i\le n\}$, является целочисленным.}
  \end{Cor}

  Аналогичные соображения позволяют также получить положительное решение проблемы~\ref{3Cycles}.

 \begin{Cor}\label{Cor5} {\it Граф ${\rm Cay}(G,S)$, где $G=A_n$~--- знакопеременная группа степени~$n$, а
  $S=R\cup R^{-1}$ для $R=\{(12i)\mid i=3,\dots,n\}$, является целочисленным.}
  \end{Cor}

  Сочетание следствий~\ref{Cor1} и~\ref{Cor3}, по-видимому, может быть полезно для нахождения новых целочисленных графов Кэли.

  \section{Обозначения}

  Мы будем использовать стандартные сведения, понятия и обозначения из теории представлений и характеров, которые могут быть найдены читателем в~\cite{Isaacs}.
 \begin{itemize}
\item[]$\mathbb{C}G$ --- комплексная групповая алгебра группы $G$.


\item[]$\Irr(G)$ --- множество неприводимых обыкновенных характеров группы $G$.

\item[]$V_\chi$ --- некоторый $\mathbb{C}G$-модуль, отвечающий (необязательно неприводимому) характеру $\chi$ группы $G$.

\item[]$\X_\chi$ --- представление группы $G$ и алгебры $\mathbb{C}G$, отвечающее характеру $\chi$.

\item[]$a_V$  для $\mathbb{C}G$-модуля $V$ и элемента $a\in \mathbb{C}G$~--- линейное преобразование пространства $V$, задаваемое правилом $v\mapsto va$.

\item[]$a_\chi$ для  элемента $a\in \mathbb{C}G$~--- линейное преобразование $a_{V_\chi}$ пространства~$V_\chi$.

\item[]$[\varphi]_B$ для линейного преобразования $\varphi$ векторного пространства $V$ c базисом $B$~--- матрица преобразования $\varphi$ в базисе~$B$.

\item[]$[\varphi]$ для линейного преобразования $\varphi$ векторного пространства $V$~--- матрица преобразования $\varphi$ в некотором базисе пространства~$V$.

\item[]$\omega_\chi$ для $\chi\in\Irr(G)$~--- гомоморфизм комплексных алгебр $Z(\mathbb{C}G)\rightarrow \mathbb{C}$, такой, что $[a_\chi]=\omega_\chi(a)I$
для любого $a\in Z(\mathbb{C}G)$ \cite[стр.~35]{Isaacs}. Здесь $I$~--- единичная матрица размера $\chi(1)\times \chi(1)$.

\item[]$\overline{S}$ для множества $S$ элементов группы $G$ --- элемент $\sum\limits_{s\in S} s\in \mathbb{C}G$.

\item[]Под {\it спектром элемента $a\in \mathbb{C}G$ на $\mathbb{C}G$-модуле $V$} будем понимать спектр линейного преобразования~$a_V$.
\item[]Любой элемент $a\in \mathbb{C}G$ мы отождествляем  с  умножением справа на $a$ элементов {регулярного модуля}. В частности, по умолчанию {\it спектр элемента $a\in \mathbb{C}G$}~--- это спектр $a$ на регулярном модуле.

 \end{itemize}

 \section{Предварительные результаты}\label{sec1}

 Напомним, что комплексное число называется {\it алгебраическим целым}, если оно является корнем некоторого многочлена с целыми коэффициентами и старшим коэффициентом 1.
 Известно, что алгебраические целые числа образуют
 подкольцо в  $\mathbb{C}$ \cite[следствие~3.5]{Isaacs}  и что число $\alpha\in \mathbb{Q}$ явяляется алгебраическим целым тогда и только тогда, когда $\alpha\in\mathbb{Z}$  \cite[лемма~3.2]{Isaacs}.

\begin{Lemma}\label{diag} Пусть $\mathfrak{X}:G\rightarrow {\rm GL}_n(\mathbb{C})$~--- представление группы $G$ с характером $\chi$ и $g\in G$~--- элемент порядка $m$. Тогда
\begin{itemize}
 \item[$(1)$] матрица $\mathfrak{X}(g)$ подобна диагональной матрице ${\rm diag}(\zeta_1,\dots,\zeta_n)$, такой, что $\zeta_i^m=1$ для всех $i=1,\dots,n$.
 \item[$(2)$] $\chi(g)=\zeta_1+\dots+\zeta_n$~--- алгебраическое целое число.
 \end{itemize}
 \end{Lemma}
 \noindent{\sc Доказательство.} Утверждение~(1) доказано в \cite[лемма~2.15]{Isaacs}. Утверждение~(2) следует из (1) и того, что алгебраические целые числа образуют подкольцо в  $\mathbb{C}$.\qed\medskip

\begin{Lemma}\label{aut} Пусть  $g\in G$~--- элемент,  порядок которого делит $m\in\mathbb{Z}$,  $\chi$~--- характер группы~$G$ и $\zeta$~--- примитивный комплексный корень степени $m$ из~$1$. Тогда:
\begin{itemize}
 \item[$(1)$]  $\chi(g)\in\mathbb{Q}(\zeta)$;
 \item[$(2)$] $\mathbb{Q}(\zeta)/\mathbb{Q}$ --- нормальное расширение;
 \item[$(3)$] для любого автоморфизма $\sigma$ поля $\mathbb{Q}(\zeta)$ существует такое $k\in \mathbb{Z}$, что ${(k,m)=1}$ и   $\zeta^\sigma=\zeta^k$;
\item[$(4)$] $\chi(g)^\sigma=\chi(g^k)$ для таких  $\sigma$ и $k$.
 \end{itemize}
 \end{Lemma}
\noindent{\sc Доказательство.} Утверждение~(1) следует из леммы~\ref{diag}~(1) и того, что $\zeta$~--- примитивный корень степени $m$ из~1.
Утверждение (2)  вытекает из того, что $\mathbb{Q}(\zeta)$~--- поле разложения
многочлена $x^m-1$.   Утверждение (3) справедливо, поскольку автоморфизм поля $\mathbb{Q}(\zeta)$ должен переводить $\zeta$ в некоторый примитивный корень степени $m$ из~$1$,
т.~е. в число вида $\zeta^k$, где  $(k,m)=1$.
Теперь если $\chi(g)=\zeta_1+\dots+\zeta_n$ , где $\zeta_i$ выбраны как в лемме~\ref{diag}~(1), то они являются корнями степени $m$ из~1,
и каждое из чисел $\zeta_i$ является степенью числа $\zeta$. Поэтому $\zeta_i^\sigma=\zeta_i^k$. Как следует из леммы~\ref{diag}, $$\chi(g)^\sigma=\zeta_1^k+\dots+\zeta_n^k=\chi(g^k).$$
Утверждение (4) доказано.\qed\medskip

\begin{Lemma}\label{omega}  Пусть $\chi\in\Irr(G)$. Тогда для любого $x\in G$ число $\omega_\chi(\overline{x^G})$ целое алгебраическое и
 $$\omega_\chi(\overline{x ^G})=\frac{\chi(x) |x^G|}{\chi(1)}.$$
\end{Lemma}
\noindent{\sc Доказательство.} См. \cite[стр.~36, в частности, теорема~(3.7)]{Isaacs}.\qed\medskip

\begin{Lemma}\label{disconnected}  Пусть $S\subseteq G$ --- симметричное множество,  $H=\langle S\rangle$ и $|G:H|=n$. Положим $\Gamma=\Cay(G,S)$ и $\Delta=\Cay(H,S)$. Обозначим через $f_G$ и $f_H$ соответственно характеристические
многочлены  матриц смежности графов $\Gamma$ и $\Delta$. Справедливы следующие утверждения.
\begin{itemize}
\item[$(1)$] Граф $\Gamma$ связен тогда и только тогда, когда $G=H$.
 \item[$(2)$] Если  $H\ne G$,  то у графа  $\Gamma$ каждая компонента связности изоморфна $\Delta$ и число комнонент равно $n$.
  \item[$(3)$] $f_G=f_H^{n}$; в частности, спектры графов  $\Gamma$ и $\Delta$ совпадают.
 \end{itemize}
\end{Lemma}
\noindent{\sc Доказательство.} Утверждение (2) доказано в \cite[лемма~1]{KonLyt}. Из утверждения (2) следуют утверждение (3) и несвязность графа $\Gamma$ при $G\ne H$ в утверждении (1).
Для завершения доказательства (1)
достаточно заметить, что любой элемент $g\in G$ представляется в виде произведения элементов множества $S$, если граф $\Gamma$ связен. Связность $\Gamma$ влечет существование пути
$$
(x_0, x_1), (x_1, x_2),\dots, (x_{m-1}, x_m)\in E(\Gamma), \text{ где } x_0=1 \text{ и } x_m=g.
$$ По определению графа Кэли для любого $i=1,\dots,m$ существует элемент $s_i\in S$, такой, что $x_i=x_{i-1}s_i$. Таким образом, $$g=x_m=x_0s_1s_2\dots s_m=s_1s_2\dots s_m.$$\qed\medskip

\section{Спектры элементов групповой алгебры}\label{sec1}

\begin{Prop}\label{AdjMatrix}{\it Пусть  $G$~--- конечная группа и $S\subseteq G$ --- такое подмножество, что $S=S^{-1}$. Тогда справедливы следующие утверждения.
\begin{itemize}
  \item[$(1)$] Матрица смежности графа $\Cay(G,S)$ совпадает с матрицей $[\overline{S}]_G$. Спектр графа $\Cay(G,S)$ совпадает со спектром элемента $\overline{S}$ на регулярном модуле.
  \item[$(2)$] Для любого элемента $a\in \mathbb{C}G$ матрица $[a]_G$ подобна блочно-диагональной матрице с блоками $[a_\chi]$ на диагонали, где $\chi\in\Irr(G)$, и каждый блок $[a_\chi]$ встречается $\chi(1)$ раз.
\item[$(3)$] Для любого элемента $a\in \mathbb{C}G$ спектр матрицы $[a]_G$ совпадает с объединением спектров матриц $[a_\chi]$ по всем $\chi\in\Irr(G)$.
\item[$(4)$] Для любого элемента $a\in \mathbb{C}G$ эквивалентны утверждения:
\begin{itemize}
\item[$(i)$] спектр элемента $a$ целочисленный $($на регулярном $\mathbb{C}G$-модуле$)$;
\item[$(ii)$] спектр элемента $a$ целочисленный на любом неприводимом $\mathbb{C}G$-модуле;
\item[$(iii)$] спектр элемента $a$ целочисленный на любом $\mathbb{C}G$-модуле.
\end{itemize}

\item[$(5)$] Если $\lambda$~--- собственное значение матрицы $[a]_G$, то кратность $\lambda$ равна $$\sum_{\chi\in \Irr(G)}\chi(1)\m_\chi(\lambda),$$ где $\m_\chi(\lambda)$~--- кратность $($возможно, нулевая$)$ $\lambda$ как собственного значения преобразования $a_\chi$.
  \item[$(6)$] Для нормального множества $S$ группы $G$ элемент $\overline{S}$ лежит в центре алгебры~$\mathbb{C}G$.

  \item[$(7)$] Для любого элемента $a\in Z(\mathbb{C}G)$ спектр матрицы $[a]_G$ равен $$\{\omega_\chi(a)\mid\chi\in\Irr(G)\}.$$

  \item[$(8)$] Пусть $x_1,\dots,x_t\in G$~--- попарно несопряженные элементы и $$S=\bigcup_{i=1}^t K_i,\quad \text{ где }\quad K_i=x_i^G.$$
  Тогда  спектр матрицы $[\overline{S}]_G$ равен $$\left\{\sum_{i=1}^t\omega_\chi(\overline K_i)\mid\chi\in\Irr(G)\right\}=\left\{\sum_{i=1}^t\frac{\chi(x_i)|K_i|}{\chi(1)}\mid\chi\in\Irr(G)\right\}.$$

\end{itemize}}
\end{Prop}
\noindent{\sc Доказательство.} Утверждение (1) легко проверить непосредственными вычислениями на основе соответствующих определений.
Утверждение (2) следует из теоремы Машке~\cite[теорема~1.9]{Isaacs} и хорошо известного разложения регулярного модуля полупростой алгебры в прямую сумму неприводимых модулей
\cite[следствие~1.17]{Isaacs}. Утверждения (3), (4) и (5) вытекают из (2). Утверждение (6) следует из~\cite[теорема~2.4]{Isaacs}. Известно, что матрица $[a_\chi]$ скалярна для любого $\chi\in\Irr(G)$ и равна $\omega_{\chi}(a)I$ (см.~\cite[стр.~35]{Isaacs}). Поэтому (7) вытекает из~(3). Поскольку отображение
$$\omega_\chi:Z(\mathbb{C}G)\rightarrow \mathbb{C}G$$
является гомоморфизмом алгебр \cite[стр.~35]{Isaacs}, в силу утверждений (1), (3), (6), (7) и леммы~\ref{omega} получаем~(8).
 \qed\medskip

\section{Доказательство теоремы~\ref{main}}\label{sec1}

\noindent{\sc Доказательство теоремы~\ref{main}.} Пусть $S$~--- нормальное эйлерово подмножество в $G$ и $$a=\overline{S}\in\mathbb{C}G.$$ В силу утверждений (1), (3), (6) и (7)  предложения~\ref{AdjMatrix} достаточно показать, что
$\omega_\chi(a)\in \mathbb{Z}$ для любого $\chi\in\Irr(G)$. Пусть $x_1,\dots,x_t$~--- представители всех классов сопряженности, объединением которых является~$S$. Тогда
\begin{equation}
 \omega_\chi(a)=\sum_{i=1}^t\frac{\chi(x_i)|x_i^G|}{\chi(1)}.\label{e1}
\end{equation}
 Из леммы~\ref{aut} следует, что $\omega_\chi(a)\in\mathbb{Q}(\zeta)$, где $\zeta$~--- примитивный корень степени $m$ из~$1$, а $m$~--- наименьшее общее кратное чисел~$|x_i|$.
Возьмем произвольный автоморфизм $\sigma$ поля $\mathbb{Q}(\zeta)$. Из леммы~\ref{aut} следует, что существует такое число $k$,
что $(m,k)=1$ и $\chi(x_i)^\sigma=\chi(x_i^k)$ для всех $i=1,\dots,t$. Поскольку множество
$S$ эйлерово, $x_i^k\in S$ для всех $i$. Ввиду нормальности множества $S$ элемент $x_i^k$ сопряжен с некоторым $x_{j}$. При этом $$\chi(x_i)^\sigma=\chi(x_i^k)=\chi(x_{j}).$$ Отображение $x\mapsto x^k$
является биекцией между классами  $x_i^G$ и $x_j^G$, поэтому
$$|x_i^G|=|x_j^G|.$$
Таким образом, автоморфизм $\sigma$ переставляет слагаемые в правой части равенства~(\ref{e1}). Тем самым  $\omega_\chi(a)^\sigma=\omega_\chi(a)$ для любого автоморфизма $\sigma$ поля $\mathbb{Q}(\zeta)$, и значит, ${\omega_\chi(a)\in \mathbb{Q}}$.
Поскольку число $\omega_\chi(a)$~--- сумма целых алгебраических чисел по лемме~\ref{omega}, оно само является алгебраическим целым и, значит, просто целым. \qed\medskip

\section{Доказательство следствий}\label{sec1}
\noindent{\sc Доказательство следствий~\ref{Cor1} и~\ref{Cor2}.}
Поскольку симметричное множество $S$, для которого $\{|s|\mid s\in S\}\subseteq\{2,3,4,6\}$, является эйлеровым, утверждение следствия~\ref{Cor1} прямо вытекает из теоремы~\ref{main}.
Следствие~\ref{Cor2} является частным случаем
следствия~\ref{Cor1}.
\qed\medskip

\noindent{\sc Доказательство следствия~\ref{Cor3}.} В групповой алгебре $\mathbb{C}G$ рассмотрим элементы $$a=\overline{S},\quad b=\overline{R}\quad \text{ и }\quad c=\overline{H\cap R}.$$
Ясно, что $a=b-c$. В силу утверждения~(6) предложения~\ref{AdjMatrix} $b\in Z(\mathbb{C}G)$, поэтому элементы $a$, $b$ и $c$ попарно коммутируют. Поскольку их матрицы симметричны, регулярный модуль $V$ групповой
алгебры обладает базисом из общих собственных векторов элементов  $a$, $b$ и $c$. Следовательно, любое собственное значение $\alpha$ элемента $a$ равно разности $\beta-\gamma$ некоторых собственных значений $\beta$ и $\gamma$ элементов
$b$ и $c$ соответственно. Спектр элемента $b$ на $V$ целый по теореме~\ref{main}. Множество $H\cap R$ является нормальным эйлеровым подмножеством в группе $H$. По  теореме~\ref{main} элемент $c$ имеет целый спектр
на регулярном модуле групповой алгебры $\mathbb{C}H$. Поскольку мы можем рассматривать любой $\mathbb{C}G$-модуль (в частности,~$V$) как $\mathbb{C}H$-модуль с согласованным действием, в силу утверждения (4)
предложения~\ref{AdjMatrix} элемент $c$ имеет на $V$ целый спектр. Значит, спектр элемента $a$ на $V$ тоже целый, и граф $\Cay(G,S)$ целочисленный по утверждению (1) предложения~\ref{AdjMatrix}.\qed\medskip

\noindent{\sc Доказательство следствия~\ref{Cor4}.} Пусть $G=S_n$, $R=\{(ij)\mid 1\le i< j\le n\}$ и $H\cong S_{n-1}$~--- множество всех подстановок в группе $G=S_n$, оставляющих символ 1 на месте. Тогда легко
заметить, что $$S=\{(1i)\mid 1< i\le n\}=R\setminus (R\cap H).$$ Как и в следствии~\ref{Cor2}, $R$~--- нормальное эйлерово подмножество в $G$. Ввиду следствия~\ref{Cor3} граф  $\Cay(G,S)$ целочисленный.
\qed\medskip

\noindent{\sc Доказательство следствия~\ref{Cor5}.} В силу леммы~\ref{disconnected} достаточно показать, что целочисленным будет граф~$\Cay(G^*,S)$, где $G^*=S_n$.
В групповой алгебре $\mathbb{C}G^*$ рассмотрим элементы
$$
a=\overline{S}=\sum_{i=3}^n\big((12i)+(21i)\big),
$$
$$
b=\sum_{1\le i<j\le n}(ij),
$$
$$
c=\sum_{3\le i<j\le n}(ij) \text{ и }
$$
$$d=(12).
$$
Согласно предложению~\ref{AdjMatrix} достаточно показать, что спектр элемента $a$ на регулярном $\mathbb{C}G^*$-модуле целый. Очевидно, что элементы $c$ и $d$ коммутируют. Элемент $b$, сумма всех транспозиций, является центральным в алгебре $\mathbb{C}G^*$, поэтому $b$, $c$ и $d$ попарно перестановочны. Непосредственно проверяется, что
$$
a=(12)\left(\sum_{i=3}^n(1i)+\sum_{i=3}^n(2i)\right)=d(b-c-d).
$$
Регулярный $\mathbb{C}G^*$-модуль обладает базисом из общих собственных векторов элементов $b$, $c$ и $d$, который будет также базисом из собственных векторов элемента $a$. При этом, если для такого вектора $v$ и чисел
$\beta,\gamma,\delta\in\mathbb{C}$ выполняются
равенства $$vb=\beta v,\quad vc=\gamma v\quad \text{ и }\quad vd=\delta v,$$ то очевидно, что $$va=\delta(\beta-\gamma-\delta)v.$$  Для завершения доказательства  достаточно показать, что спектры
элементов $b$, $c$ и $d$ целые.

Спектр $b$ целый ввиду следствия~\ref{Cor1} и утверждения (1) в предложении~\ref{AdjMatrix},
так как $b=\overline{T}$, где $T$~--- множество всех транспозиций в $G^*$. Поскольку $d$~--- элемент порядка~2, его спектр содержится среди корней многочлена $x^2-1$ и состоит из чисел~$\pm1$.  Пусть,
наконец, $H$~--- подгруппа
в группе $G^*$, изоморфная $S_{n-2}$ и состоящая из всех подстановок, оставляющих неподвижными символы 1 и~2. Тогда $H\cap T$~--- нормальное эйлерово подмножество в $H$, $c= \overline{H\cap T}$ и,
как и в доказательстве следствия~\ref{Cor3}, получаем целочисленность спектра элемента $c$. Значит, элемент $a$ также имеет целочисленный спектр. Следствие доказано.
\qed
\bigskip

Адреса авторов:
\medskip

Wenbin Guo

Department of Mathematics,

University of Science and Technology of  China,

Hefei, 230026, P. R. China

e-mail: wguo@ustc.edu.cn

\medskip

Лыткина Дарья Викторовна

Новосибирский государственный университет

ул. Пирогова, 2,

Новосибирск, 630090
Россия;

Сибирский университет телекоммуникаций и информатики,

ул. Кирова, 86,

Новосибирск, 630102 Россия.

e-mail: daria.lytkin@gmail.com

\medskip

Мазуров Виктор Данилович

Институт математики им. С.Л.Соболева СО РАН,

пр. акад. Коптюга, 4,

Новосибирск, 630090,
Россия;

Новосибирский государственный университет

ул. Пирогова, 2,

Новосибирск, 630090,
Россия.

e-mail: mazurov@math.nsc.ru

\medskip

Ревин Данила Олегович

Department of Mathematics,

University of Science and Technology of  China,

Hefei 230026, P. R. China;

Институт математики им. С.Л.Соболева СО РАН,

пр. акад. Коптюга, 4,

Новосибирск, 630090,
Россия;

Новосибирский государственный университет

ул. Пирогова, 2,

Новосибирск, 630090,
Россия.

e-mail: revin@math.nsc.ru

\end{document}